\documentclass{amsproc}

\usepackage{amsmath}
\usepackage{amssymb,amsthm}
\usepackage{stmaryrd}

\newtheorem{theorem}{Theorem}[section]

\theoremstyle{definition}
\newtheorem{definition}[theorem]{Definition}
\newtheorem{example}[theorem]{Example}

\newtheorem{corollary}[theorem]{Corollary}
\newtheorem{conj}[theorem]{Conjecture}

\theoremstyle{remark}
\newtheorem{remark}[theorem]{Remark}

\numberwithin{equation}{section}

\newcommand{\C}{\mathbb{C}}
\newcommand{\G}{\mathbb{G}}
\newcommand{\N}{\mathbb{N}}
\newcommand{\Q}{\mathbb{Q}}
\newcommand{\R}{\mathbb{R}}
\newcommand{\Z}{\mathbb{Z}}

\newcommand{\End}{\rm{End }}
\newcommand{\Endzero}{\rm{End}^0}
\newcommand{\Hg}{\rm{Hg(X) }}
\newcommand{\MT}{\rm{MT(X) }}
\newcommand{\Sp}{\rm{Sp(V,\phi) }}

\begin{document}

\title[A survey around the H., T. and M-T. conjectures for abelian varities]{A survey around the Hodge, Tate and Mumford-Tate conjectures for abelian varieties}
\author{Victoria Cantoral-Farf\'an}
\address{IMJ-PRG, UP7D B\^atiment Sophie Germain, Paris, France}
\curraddr{Institut de Math\'ematiques de Jussieu - Paris Rive Gauche (IMJ-PRG)
UP7D - Campus des Grands Moulins
B\^atiment Sophie Germain
Case 7012
75205 PARIS Cedex 13}
\email{victoria.cantoral-farfan@imj-prg.fr}
\thanks{The author was supported by a scholarship of the Consejo Nacional de Ciencia y Tecnolog\'ia, M\'exico.}

%    General info
\date{January 16, 2015 and, in revised form, ?, 2015.}

%\dedicatory{First of all I would like to thank Professor Marc Hindry for all his support and suggestions for this survey. I am immensely grateful to him for his comments on this manuscript. This research was supported supported by the Conacyt\hspace{-0.5cm}'s scholarship. I would also like to acknowledge all my colleagues from the Institut de Math\'ematiques de Jussieu – Paris Rive Gauche (IMJ-PRG) who follow closely the creation of this paper. Moreover I would like to thank the IMJ-PRG for all their support. }


\begin{abstract}
We describe Hodge, Tate and Mumford-Tate conjectures for abelian varieties. After some preliminaries on endomorphism ring, polarization and algebraic cycles, we state the three conjectures and provide a list of know results. Finally, we explain some links between these conjectures.
\end{abstract}

\maketitle
\tableofcontents

\section*{Acknowledgment}
First of all I would like to thank Professor Marc Hindry for all his support and suggestions for this survey. I am immensely grateful to him for his comments on this manuscript. This research was supported supported by the Conacyt's scholarship. I would also like to acknowledge all my colleagues from the Institut de Math\'ematiques de Jussieu – Paris Rive Gauche (IMJ-PRG) who follow closely the creation of this paper. Moreover I would like to thank the IMJ-PRG for all their support.

\section{Introduction}

The aim of this survey is to present Hodge, Tate and Mumford-Tate conjectures focusing on abelian varieties and to describe the links between them. Furthermore we will try to introduce most of the known results. Several surveys about the Hodge conjecture, for example, have been published in the past like B. Gordon's survey \cite{G} or B. van Geemen \cite{Gee}. Nevertheless few surveys about the three conjectures and their links have been written, for example W. Chi \cite{C} has presented some results in his papers concerning the three conjectures.

Actually the Hodge conjecture was introduced in 1950, its main goal is to establish a bridge between Algebraic Geometry and Differential Geometry. Hodge was inspired by Lefschetz's theorems. This conjecture is stated for algebraic varieties which are projective and smooth. Nevertheless it is easier to describe and even to establish some results about this conjecture for complex abelian varieties. One of the reasons is the simple decomposition of the singular cohomology of an abelian variety and its Hodge decomposition. Moreover we can wonder if an equivalent of this conjecture exists in Arithmetic Geometry. The answer is the Tate conjecture, stated in 1963. As for the previous conjecture, we will focus our attention on abelian varieties this time over a number field $k$. Instead of the singular cohomology we use the \'etale $\ell$-adic cohomology for abelian varieties. 
Each conjecture states, respectively, that the Hodge classes (resp. Tate classes) are $\Q$ (resp. $\Q_\ell$) linear combinations of algebraic classes. Furthermore, there exist an analogy introduced by Mumford between both conjectures. This may be stated as follows: Hodge classes (resp. Tate classes) are invariants in cohomology under the action of the Hodge group $\Hg$ (resp. Galois monodromy group $H_\ell$). After this analogy Mumford and Tate wonder if it will be relevant to compare both groups after extension of scalars to $\Q_\ell$. The answer is the Mumford-Tate conjecture, stated in 1966. 

The structure of the survey is the following. First of all we present some important notions and properties about abelian varieties: polarization, endomorphism ring. Then in section $3-5$ we introduce each conjecture in a chronological order, we present the principal ideas, introduce some general notations and describe several results. In the final section we establish some known links between those three conjectures.

%%%%%%%%%%%%%%%%%%%%%%%%%%%%%%%%%%%%%%%%%%%%%%%%%%%%%%%%%%%%%%%%%%%%%%%%%%%%%%%%%%%%%%%%%%%%%%%%%%%%%%%%%%%%%%%%

\section{Preliminaries}

The three conjectures to be discussed were formulated for algebraic varieties which are projective and smooth in general. Nevertheless, we are just going to focus our attention on abelian varieties defined over $\C$, a number field or a finite field. 

Let us remind some useful definitions concerning abelian varieties (see \cite{BL}). 

We will consider, from now on, the endomorphism ring 
$$D:=\Endzero (X)=\End (X) \otimes \Q.$$ 

Let define the polarization of an abelian variety $X$

\begin{definition}
Let $\check{X}=Pic^{\circ}(X)$ define the dual abelian variety of $X$. A polarization is a homomorphism $\lambda : X \to \check{X}$ corresponding to an ample line bundle $\mathcal{L}$ via 
$$
\lambda(x)=t_x^{*}\mathcal{L}\otimes \mathcal{L}^{-1}, \quad \forall x\in X.
$$
Any homomorphism $\alpha : X_1 \to X_2$ induce $\check{\alpha} : \check{X_2} \to \check{X_1}$.
\end{definition}

A polarization $\lambda$ of $X$ induces a positive Rosati involution $\dagger:D\to D$ defined as follows:

\begin{definition}
For every $a\in D$ we have $\check{a} \in D^{opp}$. We have the following equality $\lambda \circ a^{\dagger} = \check{a} \circ \lambda$. Moreover there exist $m\geq 1$ and $b \in \End (X)$ such that $\lambda \circ b = m \check{a} \circ \lambda$.

The Rosati involution corresponds to 
$$a^{\dagger}=\frac{1}{m} b \in D=\End(X)\otimes_{\Z}\Q.$$
\end{definition}

Consider $X$ to be a simple abelian variety, which means that every abelian subvariety of $X$ is trivial. Let us recall the Albert classification for simple abelian varieties over $\C$. 

Let us fix some notations: 

$$F=Cent(D) \quad \text{and} \quad F_0 = \{a\in F, a^{\dagger}=a\}.$$

We will now focus our attention on the following tower of extensions 
$$\Q \subset F_0 \subset F \subset D$$ 
where $[F_0 : \Q]=e_0, \quad [F : \Q]=e \quad \text{and} \quad [D:F]=d^2$.

We are going to introduce Albert's classification. We say that an abelian variety $X$ is of type (I, II, III, IV) in Albert's classification if the endomorphism ring $D=\Endzero (X)$ and the involution $\dagger$ are of one of those types. 

\begin{itemize}
\item \underline{Type I} : When $D=F=F_0$ and hence $D$ is totally real and $\dagger=id_D$
\item \underline{Type II} : When $d=2$ and $F=F_0$, and hence $D$ is an indefinite quaternion algebra ($D$ splits at all infinite places) over a totally real field $F$, which is equivalent to say that for every place $v$ at the infinity we have $inv_v (D)=0$.
Moreover there exists $a\in D-\{0\}$ with $a^2 \in F$ totally negative such that the Rosati involution is defined as follows: $\dagger : x \mapsto a^{-1} (Tr^0_{D/F} (x)-x) a$ where $Tr^0_{D/F}$ is the reduced trace map of $D/F$.
\item \underline{Type III} : When $d=2$ and $F=F_0$, and hence $D$ is a definite quaternion algebra ($D$ is inert at all infinite places) over a totally real field $F$, which is the equivalent to say that for every place $v$ at the infinity we have $inv_v (D) \neq 0$. 

Moreover $\dagger : x \mapsto Tr^0_{D/F} (x)-x$.
\item \underline{Type IV} : When $e=2e_0$ (and $d$ is arbitrary) and $F$ is a CM-field (see definition below). Moreover the restriction of $\dagger$ to $F$ is defined as the complex conjugation over $F$, $\dagger_{|F} : x \mapsto \bar{x}$.
\\
\end{itemize}

\begin{definition}
A CM-field $F$ is a totally imaginary quadratic extension of a totally real field $F'$.
\end{definition}

Let us recall that this classification concerns simple abelian varieties. Nevertheless if we consider an abelian variety $X$ isogenous to $\prod_{i}^{r} X_{i}^{n_i}$ where the $X_i$ are simple abelian varieties no isogenous then we get the following expression for the endomorphism ring of $X$:

$$
\Endzero(X)=\prod_{i}^{r} Mat(n_i \times n_i, \Endzero(X_i)).
$$

When $Y$ is an irreducible subvariety of $X$ of codimension $j$, for any cohomology theory with coefficient in $K$, it corresponds to a cohomology class $cl(Y)$ in $H^{j}(X)$, the $K$-linear span of such class form the subring of algebraic classes. The main example for us will be singular of Betti cohomology with coefficients in $\Q$ (abelian varieties defined over a subfield of $\C$) an $\ell$-adic cohomology with coefficients in $\Q_\ell$ (abelian varieties defined over a field of characteristic different from $\ell$). 

Hodge and Tate conjectures can be viewed as a description of the subring of algebraic classes while Mumford-Tate conjecture provide a bridge between these theories.

\begin{definition}
The group of algebraic classes of codimension $p$ is defined as follows:
$$C^p(X) = cl(Z^p(X)_\Q) \subset H^{2p}(X,\Q),$$ 
where $Z^p(X)_{\Q}$ is the group generated by subvarieties of $X$ of codimension $p$ and $cl: Z^p(X)_\Q \to H^{2p}(X,\Q)$ is the cycle map.
\end{definition}

%%%%%%%%%%%%%%%%%%%%%%%%%%%%%%%%%%%%%%%%%%%%%%%%%%%%%%%%%%%%%%%%%%%%%%%%%%%%%%%%%%%%%%%%%%%%%%%%%%%%%%%%%%%%%%%%

\section{The Hodge conjecture}

The Hodge conjecture ('50) for abelian varieties can be stated as follows.

\begin{conj}For every complex abelian variety $X$ the Hodge classes are $\Q$-linear combinations of algebraic classes. \end{conj}

This conjecture was already known for complex abelian varieties of dimension less or equal than $3$; it follows directly from the Lefschetz $(1,1)$-classes theorem and the Poincar\'e duality. Roughly speaking, the Lefschetz theorem states that all the Hodge classes of codimension $1$ (classes of type $(1,1)$) are algebraic classes. We will recall the definition of Hodge classes below (see definition \ref{defhc}). In general, we know that all the algebraic classes are Hodge classes but the converse is still unknown.

\subsection*{General notations}\label{gn}

Let $X_\C$ be a complex abelian variety of dimension $g$. We define $V=H_{1}(X_\C,\Q)$ as the first homology group and $\G_{m,\Q} \subset GL(V)$ as the group of homotheties. Let $\mathbb{S}$ be the restriction of scalars of $\G_{m,\C}$ from $\C$ to $\R$. We can therefore consider the next morphism:
$$
h:\mathbb{S} \to GL(V)_{\R}.
$$

First of all, we need to define the Hodge group which was introduced by Mumford; when the dimension of the abelian variety is greater than $3$ the Hodge group helps us to determine if the Hodge conjecture holds.

\begin{definition}
The \textit{Hodge group} $\Hg$ of $X$ is the smallest algebraic subgroup of $GL(V)$ defined over $\Q$ such that $h_{|U^1}$ factors through $\Hg\otimes\R$ (where $U^1\subset \mathbb{S}$ and $U^1(\R)=\{z\in\C^*, z\bar{z}=1\}$).
\end{definition}

The \textit{Mumford-Tate group} $\MT$ is the smallest algebraic subgroup of $GL(V)$ defined over $\Q$ such that $h$ factors through $\MT\otimes\R$. We can therefore relate the Hodge group and Mumford-Tate group as follows: 

\begin{equation*}\Hg=(\MT \cap SL(V))^{\circ}.\end{equation*}

We know that we can associate to every polarization of $X$ one symplectic form $\phi : V \times V \to \Q$ such that 

\begin{equation*}\Hg \subseteq Sp_{D}(V,\phi)\subseteq Sp(V,\phi),\end{equation*} where $\rm Sp_{D}(V,\phi)$ is the centralizer of the endomorphism ring $D$ in $\Sp$, the so-called \textit{Lefschetz group} of $X$.

It is well known that the cohomology group $H^n(X_\C,\Q)$ of $X$ is naturally endowed with a Hodge structure of weight $n$ and therefore we get the following decomposition:

$$H^{n}(X_\C,\Q) \otimes \C = H^{n}(X_\C,\C)=\bigoplus_{p+q=n} H^{p,q},$$
where $H^{q,p}$ is the complex conjugate of $H^{p,q}$.

It is important to notice that we have $H^n(X_\C,\C)=\bigwedge^n H^1(X_\C,\C)$ because $X$ is an abelian variety. We have also that $H^1(X_\C,\C)=T_0(X) \oplus T_0(X)^{\vee}$ which implies that 

$$
H^n(X_\C,\C)=\bigwedge^n (T_0(X) \oplus T_0(X)^{\vee})=\bigoplus_{p+q=n} (\bigwedge^p T_0(X) \wedge \bigwedge^q T_0(X)^{\vee}).
$$

One of the reasons why we focus our attention on abelian varieties is because of the simple expression of the cohomology group of $X$. The Hodge decomposition is, in general, quite difficult to express, nevertheless, for abelian varieties, this decomposition is more explicit.

We can now define the group of Hodge classes and the Hodge ring.

\begin{definition}\label{defhc}
Let $p\in \llbracket 0,g \rrbracket$, the group of Hodge classes of codimension $p$ is defined by 
$$B^{p}(X):=H^{2p}(X_\C,\Q)\cap H^{p,p} \subset H^{2p}(X_\C,\C).$$
We therefore obtain the Hodge ring as the $\Q$-graded algebra defined by
$$B^{\bullet}(X)=\oplus_{p}B^{p}(X).$$
\end{definition}

Hodge classes can also be seen as the invariants of the action of the Hodge group $\Hg$ on the cohomology ring $H^{\bullet}(X_\C,\Q)$. 

\begin{theorem}
For every $p$ we have the following equality: $$B^p (X)=H^{2p}(X_\C,\Q)^{\Hg}.$$
\end{theorem}

This is one of the reasons why Mumford introduced the Hodge group: to describe Hodge classes as invariants of a reductive algebraic group.

We define $D^{\bullet}(X)$ as the subalgebra of $B^{\bullet}(X)$ generated by divisor classes. We will say that a Hodge class $c$ is decomposable if $c\in D^{\bullet}(X)$ and that $c$ is exceptional if $c\in B^{\bullet}(X) \setminus D^{\bullet}(X)$.

Note that when exceptional classes exist, we do not know in general if they are algebraic. That is the reason why it is important to know when they exist. Nevertheless some exceptional classes have been shown to be algebraic (see \cite{Mu}). An important question is to determine when we have $$D^{\bullet}(X)\subsetneq B^{\bullet}(X).$$

One important notion is that of \textit{abelian varieties of Weil type}. Actually they allow us to determine some examples when exceptional classes may exist. We are going to clarify this point. 

\begin{definition}
Let $X$ be an abelian variety. We say that $X$ is of \textit{Weil type} (see \cite{MZ96}) when 
\begin{itemize}
\item $dim X = 2n$;
\item there exist a totally imaginary quadratic field $F$ such that $F\hookrightarrow \Endzero (X)$; 
\item and $F$ acts over the tangent space $T_0(X)$ of $X$ at $0$ with multiplicity (n,n).
\end{itemize}
\end{definition}

Let us introduce the space of Weil classes $W_F$, where $F$ is the totally imaginary quadratic field such that $F\hookrightarrow \Endzero (X)$

$$
W_F = \bigwedge_{F}^{r} V^{\vee},
$$

where $r=dim_F V^{\vee} = 2g/[F:\Q]$. \\

The $F$ vector space $W_F$ (of dimension $1$) can be identify with a subgroup of the $\Q$-vector space $H^r (X_\C,\Q)$ so the Hodge group acts over $W_F$. Therefore, because the one dimension of $W_F$ if the non trivial Weil classes are Hodge classes then all of them are either exceptional or decomposable.

The Moonen and Zarhin theorem (see \cite{MZ96}) specifies when all the Weil classes are Hodge classes. Actually the first idea came directly from Weil himself. 

\begin{theorem}
Let $X$ be a complex abelian variety and $F$ a totally imaginary quadratic field such that $F\hookrightarrow D$. If $F$ acts over the tangent space $T_0 (X)$ of $X$ at $0$ with multiplicity $(n,n)$ (where $g=2n$), then all the Weil classes are Hodge classes.
\end{theorem}

The criterion (13) of the article \cite{MZ96} allows us to determine when the Weil classes are exceptional (for instance when $X$ is of type III or IV with certain conditions).

Moreover, Moonen and Zarhin give precise conditions when the Hodge conjecture holds for simple abelian varieties over $\C$ of dimension $4$ \cite{MZ95}. Four years later they also proved under which conditions the Hodge conjecture holds for complex abelian varieties (not necessarly simple) of dimension less or equal than $5$ \cite{MZ99}. We are going to specify the previous theorems in the following section. Nevertheless, we are going to explain some useful notions.

Let $X\simeq Y^{m}$ where $Y$ is a complex abelian variety then $$\Hg = \rm Hg(Y^m) = Hg(Y).$$ Moreover, if $X\simeq \prod_i X_{i}$ where the $X_i$ are complex abelian varieties non isogenous then we have $$\Hg \subseteq \prod_i \rm Hg(X_i).$$

Notice however that there exist examples of strict inclusion, i.e. where $$\rm Hg(X_1 \times X_2) \subsetneq Hg(X_1) \times Hg(X_2).$$
Recalling the Lefschetz group introduced in section \ref{gn} we notice that we always have: $$\Hg \subset \rm Sp_D (V,\phi).$$

\subsection*{Some known results}

Let fix the condition (C) which is going to be used in the following theorems:

$$
B^{\bullet}(X^n)=D^{\bullet}(X^n) \quad \forall n \in \N^*.
$$

We refer to \cite{Mu} and \cite{Ha}.

\begin{theorem}[Murty \& Hazama '84]\label{MH}
Let $X$ be a complex abelian variety. One has 
$$
 (C) \quad \Leftrightarrow \quad 
\left\{
		\begin{array}{ll}
			X \text{ has no factors of type III}\\
			\text{in the Albert classification}\\
			\text{and}\\
			\Hg=\rm{Sp_D} (V,\phi).
    \end{array}
\right.
$$
\end{theorem}

We recall that $\rm Sp_{D}(V,\phi)$ is the centralizer of $D$ in $\Sp$.

\begin{remark}
This last theorem is quite important because when the condition (C) is verified for every $n$ then the Hodge conjecture holds. That means that Murty and Hazama proved (dependently) that this conjecture can be in many cases established looking at the endomorphism ring $D$ of the factors of $X$ and on the Hodge group. For example when the latter is the largest possible. Moreover, because of the criterion (13) of Moonen and Zarhin \cite{MZ96} exceptional classes exists when one of the simple abelian varieties is of type III in the Albert's classification.
\end{remark}

\begin{theorem}[Tankeev and Ribet '83]
Let $X$ be a simple abelian variety over $\C$ whose dimension is a prime number. Then we have ${\Hg=\rm{Sp_D} (V,\phi)}$ and the condition (C) is true.
\end{theorem}

\begin{remark}
When the dimension of a simple abelian variety is a prime number, type III in Albert's classification cannot exist, we refer to \cite{T} and \cite{R}.
\end{remark}

\begin{corollary}[Imai '76]
Let $X_1, ... , X_n$ be $n$ elliptic curves over $\C$ such that there is none isogenous to another. Let $X=\prod\limits_{i=1}^{n}X_i$, then $\Hg=\prod\limits_{i=1}^{n}\rm{Hg(X_i)}$. In particular $X$ satisfies condition (C).
\end{corollary}

We are going to focus our attention on complex abelian varieties of dimension greater than $3$. Note that the smallest abelian varieties of type III are those of dimension $4$.

\begin{theorem}[Moonen \& Zarhin '99]
Let $X$ be a complex abelian variety not necessarily simple of dimension less or equal than $5$. Then, except for some particular cases, the Hodge conjecture holds. 
\end{theorem}

The exceptional cases are explained in \cite{MZ95} and in \cite{MZ99}. We are going to give a brief description of each case. The abelian variety $X$ is of dimension $4$ in cases a-d and of dimension $5$ in e-g.

\begin{itemize}\label{cas}

\item \textbf{Case a:} $X$ is isogenous to $X_1 \times X_2$ where $X_1$ is an elliptic curve of type IV ($d=1$ and $e_0=1$), i.e. that $X_1$ is of CM-type ($\Endzero(X_1)=k$ is a imaginary quadratic field) and $X_2$ is a simple abelian variety of dimension $3$ such that $k \hookrightarrow \Endzero(X_2)$.
Given that $\Endzero(X_2)$ has an imaginary quadratic field then $X_2$ can only be of type IV (type II or III does not exist when $g=3$).

\item \textbf{Case b:} $X$ is a simple abelian variety over $\C$ of dimension 4 and type IV with $d=1$ and $e_0\in \{1,4\}$ such that $D=\Endzero(X)$ has an imaginary quadratic field $k$ which acts over $T_{X,0}$ with multiplicity  $(2,2)$. In this case $B^{\bullet}(X)$ contains exceptional Hodge classes and we do not know if they are algebraic (see \cite{MZ96}). 

\item \textbf{Case c:} $X$ is a simple abelian variety over $\C$ of dimension 4 and type III, then $D=\Endzero(X)$ is a definite quaternion algebra which has some imaginary quadratic fields ($\forall \alpha \in D-\{0\}, \quad \Q(\alpha)\subset D$ is an imaginary quadratic field).

\item \textbf{Case d:} $X$ is a simple abelian variety over $\C$ of dimension 4 such that $D=\Q$ and $\Hg \simeq SL_2 \times SL_2 \times SL_2$, then the condition (C) is not verified because $B^{\bullet}(X^2)$ contains exceptional Hodge classes, we do not know yet is they are algebraic classes. We notice here that Hodge conjecture holds for $X$ itself.

\item \textbf{Case e:} $X \sim X^{2}_{1} \times X_2$ where $X_1$ and $X_2$ are complex abelian varieties as in the case a, i.e. that $X_1$ is an elliptic curve of CM-type and $X_2$ is a simple abelian variety over $\C$ of dimension 3 of type IV such that $\Endzero(X_1)=k\hookrightarrow \Endzero(X_2)$.

\item \textbf{Case f:} $X \sim X_0 \times X_1 \times X_2$ where $X_1$ and $X_2$ are complex abelian varieties as in case a and $X_0$ is an elliptic curve not isogenous to $X_1$.

\item \textbf{Case g:} $X \sim X_1 \times X_2$ where $X_1$ is an elliptic curve of CM-type and $X_2$ is a simple abelian variety over $\C$ of dimension 4 such that $\Endzero(X_1)=k\hookrightarrow \Endzero(X_2)$ and $k$ acts over $T_{X_2,0}$ with multiplicity $(1,3)$.
\end{itemize}

Most of the cases where the Hodge conjecture is known exclude abelian varieties of type $IV$. The following result (due to Murty \cite{Mu}) gives a special condition on the dimension of the Mumford-Tate group, sufficient to imply Hodge conjecture.

\begin{theorem}\label{nondegen}
Let $X$ be a simple abelian variety over $\C$ of CM-type and dimension equal to $g$. The Mumford-Tate group $\MT$ of $X$ is, in this case, a torus of dimension less or equal than $g+1$.
If $dim \MT = g+1$ then the Hodge conjecture holds for $X$.
\end{theorem}

%%%%%%%%%%%%%%%%%%%%%%%%%%%%%%%%%%%%%%%%%%%%%%%%%%%%%%%%%%%%%%%%%%%%%%%%%%%%%%%%%%%%%%%%%%%%%%%%%%%%%%%%%%%%%%%%

\section{The Tate conjecture}

From now on we are going to consider an abelian variety $X$ over a number field $k$. The Tate conjecture ('63) for abelian varieties can be stated as follows (see \cite{T}).

\begin{conj} For every abelian variety $X$ over a number field $k$ the Tate classes are $\Q_\ell$-linear combinations of algebraic classes for every prime number $\ell$. \end{conj}

We can also consider the Tate conjecture for abelian varieties over a finite field $k$ (in which case we need to asume that $\ell \neq char(k)$) or over a finitely generated field $\Q(t_1,...,t_n)$. Tate classes are going to be defined in \ref{tc}. In general, we know that all the algebraic classes are Tate classes but the converse is still unknown. 

\subsection*{General notations}

Let $X_k$ be an abelian variety of dimension $g$ defined over a number field $k$ (or a finite field or finitely generated field), let $\bar{k}$ be the algebraic closure of $k$. We define $G_k$ as the Galois group $Gal(\bar{k}/k)$.

Let $\ell$ be a prime number and $T_\ell (X)$ be the $\ell$-adic Tate module of $X$. It is defined by 
$$T_\ell(X)=\varprojlim X[\ell^n]$$
where $X[\ell^n]$ is the Kernel of the multiplication by $\ell^n$. Then we define 

$$V_\ell=V_\ell(X):=T_\ell(X)\otimes_{\Z_\ell}\Q_\ell$$

We can therefore define the l-adic representation:

$$
\rho_\ell : G_k \to Aut_{\Q_\ell}(V_\ell)
$$

We introduced the following algebraic group $G_\ell$ defined as the Zariski closure of the image of $\rho_\ell$ (i.e. $G_\ell=\overline{\rho_\ell(G_k)}^{Zar}$) which is called the algebraic monodromy group at $\ell$. Serre proved that after some finite extension we can suppose that $G_\ell$ is a connected algebraic group (see \cite{S2}). \\
Let us introduce the identity component $G_{\ell}^{\circ}$ of $G_\ell$. Moreover we define the Galois monodromy group $H_\ell$:
$$
H_\ell =(G_\ell \cap SL(V_\ell))^{\circ}.
$$

\begin{theorem}[Serre]\label{conexe}
Let $k$ be a number field and consider the following application:
$$
\epsilon_\ell: G_k \to G_\ell / G_{\ell}^{\circ}.
$$
Then $\rm Ker(\epsilon_\ell)$ is independent of the prime number $\ell$.

Moreover there exists a finite extension $k_1$ over $k$ such that $G_{\ell,k_1}=G_{\ell,k}^{\circ}$.
\end{theorem}

Let us define $$\Endzero(X):=\End_{\bar{k}}(X)\otimes \Q.$$

We are going to introduce now the definition of Tate classes. The set $Tate^p (X)$ contains the Tate classes of codimension $p$. The first definition introduced by Tate was the following:

\begin{definition}\label{tc}
We say that $c$ is a Tate class of codimension $p$ if there exists an open subgroup $U$ of $G_k$ such that $c\in H^{2p}(X_k, \Q_\ell)^{U}$.
\end{definition}

\begin{remark}Let recall that $H^{2p}(X_k, \Q_\ell)$ is the \'etale $\ell$-adic cohomology group of the abelian variety $X$, it is defined as follows:
$$
H^{2p}(X_k, \Q_\ell)=(\varprojlim_n H^{2p}_{et}(X_k \underset{k}{\times} \bar{k}, \Z / \ell^n \Z))\otimes \Q_\ell .
$$
As the Betti cohomology of an abelian variety defined over $\C$ the \'etale $\ell$-adic cohomology admit a concrete description, this time in term of torsion points.
$$
H^1(X_k,\Q_\ell)=Hom_{\Q_\ell}(V_\ell,\Q_\ell)=V_\ell^{\vee}
$$
and
$$
H^p (X_k,\Q_\ell)=\bigwedge^p V_\ell^{\vee}.
$$
\end{remark}

Moreover we can generalize this last definition and state that 
$$
Tate^p (X):=H^{2p}(X_k, \Q_\ell)^{H_\ell}.
$$

\subsection*{Some known results}

We can present some important results for this conjecture. For instance, Faltings' theorem states as follow (see \cite{Fal}):

\begin{theorem}[Faltings '83]
Let $X$ be an abelian variety over the number field $k$. The $\ell$-adic representation $\rho_\ell$ is semisimple. That means that the identity component $G_{\ell}^{\circ}$ is reductive.
Moreover
$$
End_{G_k}(V_\ell)=End_{k}(X)\otimes_{\Z}\Q_\ell .
$$
\end{theorem}

In particular we have that for every $X_1$ and $X_2$ abelian varieties over a number field $k$ and $\ell$ a prime number, we have the following isomorphism:
$$
Hom(X_1,X_2) \otimes \Z_\ell \widetilde{\to} Hom_{\Z_\ell[G_k]}(T_\ell(X_1),T_\ell(X_2)).
$$
This last result allows us to show that the $H_\ell$ invariant of the $H^2(X_k,\Q_\ell)$ are algebraic classes. We can say that it is the $\ell$-adic analog of the Lefschetz theorem (for the Hodge conjecture).

The Tate conjecture holds for abelian varieties over a number field $k$ (or a finite field or finitely generated field) when $g$ is less or equal than $3$ or when $g$ is a prime number (see \ref{link}). Moreover, Moonen and Zarhin proved that it is also true, except some special cases, when $X$ is a simple abelian fourfold (see \cite{MZ95}). 
Most of the results that we already known follows directly from the link between the three conjectures. We refer to section \ref{link}.

%%%%%%%%%%%%%%%%%%%%%%%%%%%%%%%%%%%%%%%%%%%%%%%%%%%%%%%%%%%%%%%%%%%%%%%%%%%%%%%%%%%%%%%%%%%%%%%%%%%%%%%%%%%%%%%%

\section{The Mumford-Tate conjecture}

Hodge and Tate conjectures have an obvious similarity: they both state that cohomology classes invariant under a specific algebraic group are (linear combination of) algebraic classes. Let us remind that $X_k$ is an abelian variety over a number field $k$ or finitely generated field. Moreover we have that $V=H_1(X_\C,\Q)$, and by the comparison theorem we have that $V\otimes_\Q \Q_\ell \simeq V_\ell$, where for every prime number $\ell$, $V_\ell$ is defined through the $\ell$-adic Tate module of $X$.

Mumford-Tate conjecture establish a link between the Hodge and the Tate conjecture for abelian varieties (see \cite{Mum}), it can be stated as follows:

\begin{conj} For any prime number $\ell$ the identity component $G_{\ell}^{\circ}$ of $G_\ell$ is equal to $\MT \times_\Q \Q_\ell$.\end{conj}

This conjecture was stated in 1966 by Mumford in \cite{Mum}. Actually it was inspired by the example of abelian varieties of CM-type presented in the Shimura-Taniyama theory. We know, thanks to Mumford and Faltings results, that $\MT$ and $G_\ell$ are reductive groups.\\
Moreover $Z^{\circ}(\MT_{\Q_\ell})=Z^{\circ}(G_\ell)$ for every prime number $\ell$ (see \cite{V}), so in order to establish Mumford-Tate conjecture we just need to look at the semi-simple part of each of these groups. We notice as well that when $X$ is of CM-type, the semi-simple part is trivial, so we get the Mumford-Tate conjecture directly. \\
The Hodge group, as we already said, was introduced by Mumford in the same article. Actually it was introduced to show more clearly the link between Hodge and Tate conjectures and then to state the Mumford-Tate conjecture.

Following the same ideas, Serre also proved that the rank of the semi-simple part of $G_\ell$ does not depend of the prime number $\ell$. 

Deligne (Piatetskii-Shapiro, Borovoi) proved that we have for every prime number $\ell$ the first inclusion: $$G_{\ell}^{\circ} \subseteq \MT \times_\Q \Q_\ell .$$

Moreover Serre also proved that if there exists one prime number $\ell$ such that $rk_{ss}(G_\ell)\geq rk_{ss}(\MT)$ then, because of the Deligne inclusion and the fact that they are both reductive groups, we obtain the equality $G_{\ell}^{\circ}=\MT \times_\Q \Q_\ell$ for every prime number $\ell$. The fact that the rank is independent to $\ell$ allows us to prove that if this equality holds for just one $\ell$ it holds for every $\ell$. This last result is due to Larsen and Pink in 1995.

\subsection*{General notations}

Let $X_k$ be a simple abelian variety of dimension $g$ over a number field $k$ (or a finite field or finitely generated field). Let us remind that $[F:\Q]=e$ and $[D:F]=d^2$ where $F=Cent(D)$ and $D$ is the endomorphism ring. We define the relative dimension of $X$ as follow:

\begin{definition}
$$
\rm reldim(X)=\left\{
		\begin{array}{ll}
			\frac{g}{e} \text{if $X$ is of type I}\\
			\frac{g}{2e} \text{if $X$ is of type II or III}\\
			\frac{2g}{de} \text{if $X$ is of type IV}
    \end{array}
\right.
$$
\end{definition}

\subsection*{Some known results}

We notice, like for the Hodge and Tate conjecture, that Mumford-Tate conjecture holds when $X$ is a simple abelian variety dimension less or equal than $3$. Furthermore, Mumford-Tate conjecture holds for simple abelian varieties of prime dimension, see \cite{C1}.

We are going to present in a chronological order some important results concerning the Mumford-Tate conjecture. 

First of all, Serre proved in 1972 that Mumford-Tate conjecture holds for elliptic curves. Then, before Faltings' theorem, Ribet presented a particular case where the Mumford-Tate conjecture holds. This one occures when $D$ is a totally real field of dimension $g$ and there exists a place of bad semi-stable reduction.

After Faltings' result, Serre, Chi, Pink and Banaszak, Gajda and Kraso\'n presented new cases where the Mumford-Tate conjecture holds for abelian varieties of type I and II. We are going to introduce each result.
\\

The following theorem of Serre was presented during his course at the Coll\`ege de France in 1984-1985 (see \cite{S}). 

\begin{theorem}[Serre]
If we suppose that $\Endzero(X)=\Q$ and $\rm dim(X)$ is an odd integer then the Mumford-Tate conjecture holds. 
\end{theorem}
\begin{remark}
The proof of this theorem can be found in \cite{C} section 6. Moreover Serre proved with this same hypothesis that Hodge and Tate conjecture hold too. 
Let also remark that when $dim X=2$ and $\Endzero(X)=\Q$, Mumford-Tate holds. Actually in this case $\Hg = Sp_4$ (because $\Hg$ is an irreducible reductive group inside $\rm Sp_4$) and so does $G_{\ell}^{\circ}$. Serre also proved the same for $g=6$.
\end{remark}

Chi improves Serre results using the same technique, he introduced the next theorem, which covers many cases when $X$ is of type I or II.

\begin{theorem}[Chi '92]
The Mumford-Tate conjecture hold for this four cases:
\begin{itemize}
\item When $X$ is an abelian variety of type I,
\begin{itemize}
\item $dim X = d$ with d=2 or an odd number so that $\Endzero(X)=\Q$;
\item $dim X = 2d$ with d an odd number so that $\Endzero(X)=E$ is a real quadratic field.
\end{itemize}
\item When $X$ is an abelian variety of type II,
\begin{itemize}
\item $dim X = 2d$ with d=2 or an odd number so that $\Endzero(X)=D$ is an indefinite quaternion algebra over $\Q$;
\item $dim X = 4d$ with d an odd number so that $\Endzero(X)=D$ is an indefinite quaternion algebra over a real quadratic field.
\end{itemize}
\end{itemize}
\end{theorem}

Let us define the following set introduced by Pink where $k$ is an odd integer. We have
$$
\mathcal{S}:=\left\{g\geq 1; \exists k\geq 3, \exists a\geq 1, g=2^{k-1}a^k \right\} \bigcup \left\{g\geq 1; \exists k\geq 3, 2g=\binom{2k}{k}  \right\}.
$$

\begin{theorem}[Pink '98]
Assume that $\Endzero(X)=\Q$ and that $g \notin \mathcal{S}$ then for every prime number $\ell$ we have $G_{\ell}^{\circ} = \MT \times_\Q \Q_\ell$, in particular Mumford-Tate conjecture holds for $X$.
\end{theorem}
\begin{remark}
This last result is more related to representation theory and minuscule weights introduced by Serre, the proof is presented in \cite{P}.
\end{remark}

Few years later in 2003, Banaszak, Gajda and Kraso\'n proved Mumford-Tate conjecture for abelian varieties of type I or II see \cite{BGK}. They used Pink results to generalized Chi's theorem.

\begin{theorem}[Banaszak, Gajda, Kraso\'n]\label{bgk}
Let $X$ be an abelian variety over a number field $k$ of type I or II in Albert's classification and odd relative dimension (or equal to $2$). Then $\MT = \rm GSp(V,\phi)$ and $G_{\ell}^{\circ}=\rm GSp(V,\phi) \otimes \Q_\ell$ hence the Mumford-Tate conjecture holds for $X$.
\end{theorem}
\begin{remark}
Moreover, Hodge conjecture holds for abelian varieties with the same hypothesis as before. That means that Tate conjecture also holds for $X$ (see \ref{link}).
\end{remark}

There is another link with this conjecture that does not only concern the endomorphism ring of $X$ or its dimension. Let us present the following result:

\begin{theorem}[Hall '08]
Let $X$ be an abelian variety over a number field $k$ such that $\Endzero(X)=\Q$ and assume that the N\'eron model of $X$ over $\mathcal{O}_k$ has a semi stable fiber with toric dimension equal to $1$. Under those conditions we have that $G_\ell = GSp_{\Q_\ell}$ and Mumford-Tate conjecture holds for $X$.
\end{theorem} 
\begin{remark}
It is interesting to notice that Hall's theorem is one of the first theorems that provides information about the algebraic monodromy group $G_\ell$ and then allows us to obtain information about the Mumford-Tate group $\MT$.
\end{remark}

Moreover it is interesting to notice that M. Hindry and N. Ratazzi generalized Hall's result to abelian varieties of type I (resp. type II) with one semi stable fiber with toric dimension equals to $e_0$ (resp. $2e_0$), we refer to \cite{HR}.\\

More results have been proved since then. For instance we have two more result from last year. \\

Zhao has announced a proof of the Mumford-Tate conjecture for simple abelian fourfolds over a number field $k$ see \cite{Z14}. Actually the difficult case remaining was to prove that Mumford-Tate hold when the abelian variety $X$ is of type I in the Albert classification and $\Endzero(X)=\Q$. Because Moonen and Zarhin had already proved the other cases (see \cite{MZ95}): the idea is because of the dimension of $X$ the groups $G_\ell$ and $\MT$ do not have several option, they must be reductive groups inside a symplectic group. We are going to present an idea of the difficult part for the Mumford-Tate conjecture.

Let assume that $X$ is a simple abelian fourfold over a number field $k$ of type I such that $\End_{\bar{k}}(X)=\Z$. Under this conditions we have several options for the Hodge group:
\begin{enumerate}
\item $\Hg = Sp_{8,\Q}$;
\item $\Hg = (SL_2)_{\Q}^3$ (a twist of).
\end{enumerate}
We have the same options tensorized by $\Q_\ell$ for $H_\ell:=(G_\ell \cap SL_8)^{\circ}$ for every prime number $\ell$. Let remind that it is equivalent to the Mumford-Tate conjecture using the Hodge group. 
Let recall the following inclusions:
$$
H_\ell \subseteq \Hg \otimes \Q_\ell \subseteq \rm Sp_{8,\Q_\ell}. 
$$
We have then different possible combinations:
\begin{itemize}
\item if $H_\ell = Sp_{8,\Q_\ell}$ then $\Hg=Sp_{8,\Q}$;
\item if $\Hg = (SL_2)_{\Q}^3$ then $H_\ell = (SL_2)_{\Q_\ell}^3$.
\end{itemize}
Those combinations implies Mumford-Tate conjecture.
Nevertheless the difficult part is the following. Let suppose that $\Hg=Sp_{8,\Q}$, the question is to know if we have necessarily $H_\ell = Sp_{8,\Q_\ell}$ i.e. can we exclude the case where $H_\ell = (SL_2)_{\Q_\ell}^3$? The answer is positive as Zhao announced in \cite{Z14}. The idea is to use Shimura varieties of $\Hg$ and $H_\ell$.\\

The classic idea of Shimura varieties is to classify abelian varieties through their endomorphism rings. Consequently, a Shimura variety $\mathcal{A}_{g} (\C)$ is a family of complex abelian varieties of dimension $g$ which are characterized by their endomorphism ring. Mumford's idea (see \cite{Mum}) was to see Shimura varieties as a family of complex abelian varieties characterized by their Hodge group.
Zhao idea was to focus his attention on the Shimura varieties $\mathcal{A}_{4} (\C)$ and more particularly the family of abelian varieties (up to isogeny) such that their Hodge group is contained in $SL_{2}^{3}$. Furthermore, thanks to the previous inclusion $\rm H_\ell \subseteq Hg \otimes \Q_\ell$ we can obtain Mumford-Tate conjecture with those Shimura varieties.\\

Moreover, thanks to the work of Moonen and Zarhin in \cite{MZ99} Lombardo proved in \cite{L} that Mumford-Tate conjecture holds for abelian varieties not necessarily simple of dimension less or equal to $5$.

%%%%%%%%%%%%%%%%%%%%%%%%%%%%%%%%%%%%%%%%%%%%%%%%%%%%%%%%%%%%%%%%%%%%%%%%%%%%%%%%%%%%%%%%%%%%%%%%%%%%%%%%%%%%%%%%

\section{Links between these three conjectures \label{link}}

Let $X_k$ be an abelian variety over a number field $k$ of dimension $g$. Let us fix an embedding $\sigma: k \hookrightarrow \C$ such that $X_\C = X_k \underset{k,\sigma}{\times} \C$ then we define $V=H_1(X_\C,\Q)$. We remind that $D=\Endzero(X)$  (resp. $D_{\Q_\ell}$) can be seen as a subalgebra of $\End_\Q (V)$ (resp. $\End_{\Q_\ell} (V)$) as follow:

$$
D=\End_\Q (V)^{\rm Hg(X_\C)} \quad D_{\Q_\ell}=\End_{\Q_\ell} (V)^{G_\ell^{\circ}} 
$$ 

Let us recall the Hodge and Tate conjectures, they both state that Hodge classes or Tate classes are $\Q$ or $\Q_\ell$ (for every prime number $\ell$) linear combinaison of algebraic classes. Moreover we can present the following analogy:

$$
B^p (X)=H^{2p}(X_\C,\Q)^{\Hg}
$$

$$
Tate^p (X)=H^{2p}(X_k,\Q_\ell)^{H_\ell}
$$

It is interesting to recall that Tate stated his conjecture in 1963 while Mumford introduced the Hodge group in 1966, the same year when he stated Mumford-Tate conjecture. Let us notice that Mumford was inspired by the idea that Tate classes are classes of $H^{2p}(X_k,\Q_\ell)$ which are invariant under the action of the l-adic monodromy group $H_\ell$. 

We are going to describe links between the Mumford-Tate conjecture (MT), the Hodge conjecture (H) and the Tate conjecture (T).

We have the following equivalence:
\begin{equation}(MT)+(H) \Leftrightarrow (T).\label{equi}\end{equation}

Let us suppose that Mumford-Tate and Hodge conjectures hold for an abelian variety $X$. We have then that $G_\ell=\MT \otimes \Q_\ell$ or equivalently $H_\ell=\Hg \otimes \Q_\ell$. Moreover for every $p\in \llbracket 0,g \rrbracket$ we have that the classes in $H^{2p}(X_\C,\Q)^{\Hg}$ are algebraic classes. That means, because of the Mumford-Tate conjecture, that $H^{2p}(X_k,\Q_\ell)^{H_\ell}$ are also algebraic classes for every prime number $\ell$. This last property implies Tate conjecture. We have then that $(MT)+(H) \Rightarrow (T)$.

\begin{example} Let us present several examples where this implication is fruitful:
\begin{enumerate}
\item We know that Mumford-Tate conjecture holds for simple abelian varieties of prime dimension over a number field $k$ (see \cite{C1}) and that the Hodge conjecture also holds for simple abelian varieties over $\C$ of prime dimension (Tankeev). That means that Tate conjecture holds for this abelian varieties.
\item Furthermore, we know that Mumford-Tate conjecture holds for simple abelian varieties of CM-type (Shimura-Taniyama theory). Thanks to the theorem \ref{nondegen} Hodge conjecture holds for simple abelian variety over $\C$ of CM-type of dimension equal to $g$ and such that the Mumford-Tate group $\MT$ of $X$ is a torus of dimension $g+1$. Therefore Tate conjecture also holds for this simple abelian varieties.
\item Similarly, thanks to the theorem \ref{bgk}, we know that if $X$ is an abelian variety over a number field $k$ of type I or II and with odd relative dimension we have that $\MT = \rm GSp(V,\phi)$ and $G_{\ell}^{\circ}=\rm GSp(V,\phi) \otimes \Q_\ell$. That means that the Mumford-Tate conjecture holds for this kind of abelian varieties. Consequently we have that $\Hg:=(\MT \cap SL(V))^{\circ}=Sp(V,\phi)$. Using the theorem of Murty and Hazama \ref{MH} we know that the Hodge conjecture also holds for this abelian varieties. That implies that Tate conjecture also holds for abelian varieties over a number field $k$ of type I or II and with odd relative dimension.
\end{enumerate}
\end{example}

Piatetskii-Shapiro had proven that $(T)\Rightarrow (H)$. The ideas behind this are clearly explained in a AIM talk of Milne (\cite{M}) using absolute Hodge classes. Moreover one can also prove that $(H) + (T) \Rightarrow (MT)$ using motives.

Finally using this last three implications we can deduce the equivalence \ref{equi}. That means that we can deduce Tate conjecture for every abelian variety which satisfies the Mumford-Tate and the Hodge conjecture.

\bibliographystyle{amsalpha}

\end{document}